\documentclass[12pt]{article}
\usepackage[a4paper,left=1in,top=1in,width=6.3in,height=9in]{geometry}
\usepackage{natbib}

\usepackage{type1cm}        
\usepackage{graphicx}        

\usepackage[width=.75\textwidth]{caption}

\usepackage{amsbsy}
\usepackage{amsfonts}
\usepackage{amsmath}

\usepackage{bm}
\newcommand{\nil}[1]{}
\newcommand{\diag}{\text{diag}}

\newcommand{\ACG}{\text{ACG}}
\newcommand{\WC}{\text{WC}}
\newcommand{\Arg}{\text{Arg}}
\newcommand{\vecc}{\text{vec}}
\newcommand{\SSigma}{\bm \Sigma}
\newcommand{\OOmega}{\bm \Omega}
\renewcommand{\phi}{\varphi}
\newcommand{\Mobius}{M\"obius }

\begin{document}
\title{Directional distributions and the half-angle principle}
\author{John T. Kent\\ University of Leeds, Leeds LS2 9JT, UK\\
  \texttt
  {j.t.kent@leeds.ac.uk}
}
%
\maketitle

\abstract{Angle halving, or alternatively the reverse operation of
  angle doubling, is a useful tool when studying directional
  distributions.  It is especially useful on the circle where, in
  particular, it yields an identification between the wrapped Cauchy
  distribution and the angular central Gaussian distributions, as well
  as a matching of their parameterizations.  The operation of angle
  halving can be extended to higher dimensions, but its effect on
  distributions is more complicated than on the circle.  In all
  dimensions angle halving provides a simple way to interpret 
  stereographic projection from the sphere to Euclidean space.}

\vspace{5mm}

\section{Introduction}
\label{sec:intro}
The wrapped Cauchy (WC) distribution on the circle is a remarkable
distribution that appears in a wide variety of seemingly unrelated
settings in probability and statistics.  
The ACG distribution is another important distribution in directional
statistics.  It was used by \citet{Tyler87a,Tyler87b} to 
construct and study a  robust estimator of a covariance matrix, or more
generally a scatter matrix, for $q$-dimensional multivariate data.

As noted in \citet{kc88c}, the ACG distribution in $q=2$ dimensions
(i.e. on the circle) can be identified with the WC distribution after
angle doubling.  Equivalently, WC distribution can be identified with
the ACG distribution after angle halving.  Hence algorithms to
estimate the parameters of one distribution can be used with little
change to estimate the parameters of the other distribution.  Several
algorithms to compute the maximum likelihood estimates based on the EM
algorithm have been explored in \citet{kc88c} and \citet{kc94a}.  See
also \citet{kc95b} for further discussion.

The current paper extends the analysis as follows:
\begin{itemize}
\item to use angle halving on the circle to recast the \Mobius
  transformation in terms of a rescaled linear transformation of the
  plane, a result which additionally allows us to match the
  parameterizations of the WC and ACG distributions;

\item to extend angle halving to higher dimensions and to show the
  connection between gnomonic projection and stereographic projection;

\item to note that the ACG distribution under gnomonic projection maps
  to a multivariate Cauchy distribution; and to contrast it with the
  spherical Cauchy distribution of \cite{Kato-McC20}, which under
  stereographic projection maps to a multivariate $t$-distribution;

  \item to summarize some further properties of the WC distribution.
  \end{itemize}

To set the scene for the main investigation of the paper, recall some
basic properties of the WC and ACG distributions on the circle $S_1$.
The WC distribution, written WC$(\lambda)$,  has probability density
function (p.d.f.)
\begin{equation}
  \label{eq:wc}
  f_\WC(\theta; \lambda) = (2\pi)^{-1} \frac{1-\lambda^2}
  {1+\lambda^2  - 2 \lambda \cos \theta},
  \quad \theta \in S_1.
\end{equation}
Here $0 \leq |\lambda|< 1$ is a concentration parameter. The distribution has
been centered to have its mode at $\theta=0$ if $\lambda>0$ and $\theta=\pi$
if $\lambda<0$; it reduces to the uniform distribution if $\lambda=0$.

The ACG distribution on $S_1$, written ACG$(b)$, has probability
density function (p.d.f.)
\begin{align}
  f_\ACG(\phi; b)
  &= (2 \pi)^{-1} b/ \{ b^2 \cos^2 \phi + \sin^2 \phi\} \notag \\
  &= \pi^{-1} b/ \{ b^2 (1 + \cos 2\phi) + (1 - \cos 2\phi)\} \notag\\
  &= \pi^{-1} b / \{(1+b^2) - (1-b^2) \cos 2\phi\}, \ \phi \in S_1. \label{eq:acg2}
\end{align}
Here $0<b<\infty$ is a concentration parameter.  The density is
antipodally symmetric, $f(\theta) = f(\theta+\pi)$.  The distribution
has been centered to have its modes at $\theta=0, \pi$ if $b<1$ and
$\theta=\pm\pi/2$ if $b>1$; it reduces to the uniform distribution if
$b=1$.

If
\begin{equation}
  \label{eq:lambda-to-b}
  b=(1-\lambda)/(1+\lambda),
\end{equation}
it can be checked that (\ref{eq:acg2}) is the same as (\ref{eq:wc})
under the angle doubling relation $\theta=2\phi$.  That is, if $\Phi$
is a random angle following the ACG$(b)$ distribution and
(\ref{eq:lambda-to-b}) holds, then $\Theta=2\Phi$ is a random angle
following the WC$(\lambda)$ distribution.  The relation (\ref{eq:lambda-to-b})
between $b$ and $\lambda$ will be assumed throughout the paper.

The paper is organized as follows.  Basic transformations of the
circle are defined and examined in Section \ref{sec:circle-ops}.
These transformations are used in Section \ref{sec:circle-trans} to
obtain the ACG and WC distributions on the circle as transformations
of the uniform distribution.  The basic transformations are extended
to the sphere in Section \ref{sec:sphere-ops} and interpreted through
two projections in Section \ref{sec:sphere-proj}.  The transformations
are used to obtain the ACG distribution on the sphere (Section
\ref{sec:acg}) and a spherical analog of the WC distribution (Section
\ref{sphere:sc}) as transformations of the uniform distribution.
Finally, Section \ref{sec:wc-params} summarizes some further
derivations and motivations for the WC distribution on the circle.

For some standard background on directional distributions, see, e.g.,
\citet{Mardia-Jupp00} and \citet{Chikuse03}.  For basic results from
multivariate analysis, see, e.g., \citet{ka79}.  A fundamental
reference is \citet{McC96}, which goes further than the current paper
in exploring how the family of WC distributions is closed under the group
of \Mobius transformations on the unit circle.

\section{Basic operations on the circle}
\label{sec:circle-ops}
A point on the circle can be written as an angle $\phi$, where without loss of
generality, $\phi \in (-\pi,\pi]$.  The point can also be expressed as
a unit vector
\begin{equation}
  \label{eq:circ-rep}
\bm x  = (x_1,x_2)^T = (\cos \phi, \sin \phi)^T =
\pm (1,r)^T / \sqrt{1+r^2}, \quad r = \tan \phi,
\end{equation}
or as a complex number $x_1 + ix_2 = C(\bm x)$.
It is convenient to denote the mappings between vector
and angular representations by
\begin{equation}
  \label{eq:circ-vec-ang}
  \phi = \Arg(\bm x), \quad  \bm x=\vecc(\phi).
  \end{equation}

For later use note that the derivatives of the mappings between $\phi$
and $r=\tan \phi$ are given by
\begin{equation}
\label{eq:circ-der-r}
dr/d\phi = \sec^2 \phi = 1/ \cos^2 \phi = 1+r^2, \quad
d\phi / dr = 1/(1+r^2).
\end{equation}

Another important representation of an angle, where this time the
angle is denoted $\theta$, is in terms of the tangent of the
half-angle, $s = \tan(\theta/2)$.  Square both sides and use the
double angle formulas to get
\begin{equation}
  \label{eq:circ-s}
s^2 = \tan^2 (\theta/2) = \frac{\sin^2(\theta/2)}{\cos^2(\theta/2)}
= \frac{1-\cos \theta}{1+\cos \theta},
\end{equation}
which can be inverted to give
$$
\cos \theta = \frac{1-s^2}{1+s^2},
$$
so that  $1+\cos \theta = 2/(1+s^2)$.

Throughout the paper we assume that $\theta$ and $\phi$ are related by
the double angle condition, $\theta=2\phi$, so that $r=s$.  However,
it is helpful to use both notations $r$ and $s$ to emphasize that $r$
is obtained from $\phi$ and $s$ is obtained from $\theta$.

Three important mappings from $S_1$ to itself are as follows.
  \begin{itemize}
  \item[(a)] \emph{Squaring}, denoted $S(\bm x)$.  In vector form the
    transformation is defined by
    \begin{equation}
      \label{eq:squaring}
    S(\bm x) = (x_1^2 - x_2^2, 2 x_1 x_2)^T, \quad \bm x \in S_1.
    \end{equation}
    If $\bm y = S(\bm x)$, then in complex arithmetic
    $y_1+i y_2 = (x_1 + i x_2)^2$.  Further, if $\phi=\Arg(\bm x)$ and
    $\theta = \Arg(\bm y)$ are the two points in angular coordinates,
    then $\theta = 2\phi$.  Hence squaring is a two-to-one mapping of
    $S_1$ to itself.

  \item[(b)] The \emph{rescaled diagonal linear transformation}, denoted
    $L(\bm x; b)$, where $ b >0$ is a 
    scaling constant.  In vector form the
    transformation is defined by
    \begin{equation}
      \label{eq:BLT}
    L(\bm x; b) = (x_1, b x_2)^T / \sqrt{x_1^2 + b^2 x_2^2}.
    \end{equation}
    That is, the second component of $\bm x$ is scaled by a factor $b$, and the
    resulting vector is rescaled to be a unit vector.  
    The rescaled diagonal linear transformation can also be described as
    follows.  If $\bm z = L(\bm x; b)$ then
    \begin{equation}
      \label{eq:BLT-tan}
    \tan \Arg(\bm z) = b\, \tan \Arg(\bm x).
    \end{equation}
    
  \item[(c)] The \emph{diagonal \Mobius transformation}, denoted
    $M(\bm y; \lambda)$. In vector form the transformation is defined
    for $\lambda>0$ by
    \begin{equation}
      \label{eq:Mobius}
      M(\bm y; \lambda) = (2 \lambda + (1+\lambda^2) y_1, (1-\lambda^2) y_2)^T /
      (1+\lambda^2 + 2 \lambda y_1), \quad \bm y \in S_1.
    \end{equation}
    If $\bm w = M(\bm y; \lambda)$ where $\Arg(\bm y)=\theta$ and $\Arg(\bm w)=\eta$,
    then $\theta$ and $\eta$ are related by
    \begin{equation}
      \label{eq:Mobius-tan}
    \tan \eta/2 = b \tan \theta/2,
  \end{equation}
  where $b$ and $\lambda$ are related by (\ref{eq:lambda-to-b}).  That
  is, the \Mobius transformation is the same as the rescaled diagonal
  linear transformation after the angles $\theta$ and $\eta$ are
  divided by 2.  The \Mobius transformation is most commonly defined
  using complex arithmetic,
    \begin{equation}
      \label{eq:Mob-complex}
      C(M(\bm y; \lambda))=\frac{y_1 + i y_2 +\lambda}
      {\lambda (y_1 + iy_2) + 1},  \quad \bm y \in S_1
    \end{equation}
    where for our purposes here, $0< \lambda< 1$ is restricted to
    being real.
    \end{itemize}

  These transformations can be combined to give the following result,
  which it is helpful to call the  \emph{fundamental diagonal \Mobius identity}:
  \begin{equation}
      \label{eq:ID}
      M(S(\bm x); \lambda) = S(L(\bm x; b)), \quad \bm x \in S_1,
    \end{equation}
    where  $b$ and $\lambda$ are related by (\ref{eq:lambda-to-b}).
  That is, a rescaled diagonal linear transformation followed by squaring is the
  same as squaring followed by a diagonal \Mobius transformation.

  The identity in (\ref{eq:ID}) has been stated for diagonal case.
  However, it is possible to construct a more general version by
  allowing rotations before and after the relevant transformation.
  Let
  \begin{equation}
    \label{eq:rotation}
  \bm R_\alpha = \begin{bmatrix} \cos \alpha & -\sin \alpha \\
    \sin \alpha & \cos \alpha \end{bmatrix}
  \end{equation}
  denote a $2 \times 2$ rotation matrix by an angle $\alpha$.  Also,
  recall that any $2 \times 2$ matrix $\bm B$ with positive determinant
  can be writing using the singular value decomposition as
  $$
  \bm B = c \bm R_\alpha \diag(1,b) \bm R^T_\beta
  $$
  where $c>0$ and $b>0$.  Note that
  $\bm R_\alpha^T \vecc(\phi) = \vecc(\phi-\alpha)$ and
  $S(\bm R^T_\alpha \bm x) = \bm R^{2T}_\alpha S(\bm x) = \vecc(2(\theta-\alpha))$.

  Define more general versions of the rescaled diagonal linear and \Mobius 
  transformations by
  \begin{align}
    &L(\bm x; \bm B) = \bm B \bm x / ||\bm B \bm x|| =
                      \bm R_\alpha  L(\bm R^T_\beta\bm x; b), \notag\\
    &M(\bm x; \lambda,\exp(2i\alpha),\exp(2i\beta))) =
    \bm R_\alpha^2 M(\bm R_\beta^{2T}\bm x; \lambda), \label{eq:trans-general}
  \end{align}
  where $||\bm x||^2 = \bm x^T \bm x$.  In complex notation,
  the \Mobius transformation becomes
  $$
  M(\bm y; \lambda, \exp(2i\alpha),\exp(2i\beta))) =
  \exp(2i(\alpha-\beta)) \frac{y_1+iy_2 + \lambda \exp(2i\beta)}
  {\lambda \exp(-2i\beta)(y_1+iy_2) + 1}.
  $$
  Note the $L$ now depends on the matrix $\bm B$ and $M$ now depends
  on a real number and two complex numbers.  The more general version
  of the fundamental \Mobius identity becomes
  \begin{equation}
    \label{eq:ID-general}
  M(S(\bm x); \lambda, \exp(2i\alpha),\exp(2i\beta)) = S(L(\bm x; \bm B)).
\end{equation}

\section{Transformations of distributions on the circle}
\label{sec:circle-trans}
Let $\Phi^*$ follow a uniform distribution on the circle, with density
$f(\phi^*) = 1/(2\pi), \ -\pi < \phi^* < \pi$.  Let $R^*=\tan \Phi^*$ and
$\bm X^* = \vecc(\Phi^*)$ denote the corresponding tangent of the angle
and the Euclidean coordinates.  Consider the rescaled diagonal linear
transformation $\bm X = L(\bm X^*; b)$, where b>0, and let
$\Phi= \Arg(\bm X)$ and $R = \tan(\Phi)$ denote the corresponding
angular and tangent values.

The inverse transformation between $\bm X$ and $\bm X^*$ is
$\bm X^* = L(\bm X; 1/b)$.  Then the p.d.f. of $\Phi$ is given by
\begin{align}
  \frac{1}{2\pi} \frac{d\phi^*}{d \phi}
  &= \frac{1}{2\pi} 
    \frac{d\phi^*}{dr^*} \frac{dr^*}{dr} \frac{dr}{d\phi} \notag \\
  &=\frac{1}{2\pi}\, \frac{1}{1+r^{*2}} b^{-1} (1+r^2) \notag \\
  &= \frac{1}{2\pi b}\, \frac{\cos^2 \phi}{\cos^2 \phi + b^{-2} \sin^2 \phi} \,
    \frac{1}{\cos^2 \phi}   \notag \\
  &= \frac{b}{2\pi}\, \frac{1}{b^2 \cos^2 \phi + \sin^2 \phi} = f_\ACG(\phi;b),
    \label{eq:pdf-circ-acg}
\end{align}
where we have used the fact that
$r^{*2} = b^{-2} r^2 = b^{-2}\sin^2 \phi/\cos^2 \phi$, and
$1/(1+r^2) = \cos^2 \phi$.  In other words $\Phi$ follows the
ACG$(b)$ distribution.

If $\Phi^*$ follows a uniform distribution, then so does
$\Theta^* = 2\Phi^*$.  Hence
$$
\Theta = \Arg(M(\Theta^*, \lambda)) = 2 \Phi = 2 \Arg(L(\Phi^*, b))
$$
has p.d.f. (\ref{eq:pdf-circ-acg}) as a function of $\phi$ (the factor
1/2 from the Jacobian $d \phi^* / d \theta^*$ cancels the factor 2 which
arises since the mapping from $\phi^*$ to $\theta^*$ is two-to-one).
After writing the p.d.f. in terms of $\theta$, the wrapped Cauchy
density $f_\WC(\theta;\lambda)$ in (\ref{eq:wc}) is obtained, where
$\lambda$ is related to $b$ by (\ref{eq:lambda-to-b}).

In particular, if $0<\lambda<1$, i.e. $0<b<1$, the diagonal \Mobius
mapping $\bm Y = M(\bm Y^*, \lambda)$ pulls probability mass towards the
direction $\theta=0$; similarly the rescaled diagonal linear mapping
$\bm X = L(\bm X^*; b)$ pulls probability mass towards the directions $\phi=0$
and $\pi$.  Hence the WC distribution for $\bm Y$ has
a mode in the zero direction and the ACG distribution for $\bm X$ has
its modes in the directions $0$ and $\pi$.

In summary, both the ACG and WC distributions can be obtained from
suitable transformations of the uniform distribution.  For simplicity,
attention has been focused on the centered distributions in this
section, but rotations of the modal direction can be easily included.

\section{Basic operations on the sphere}
\label{sec:sphere-ops}
In higher dimensions, more notation is needed.  For $q \geq 2$, let
$S_{q-1}= \{\bm x \in \mathbb{R}^q : \bm x^T \bm x = 1\}$ denote the
unit sphere in $\mathbb{R}^q$ with surface area
\begin{equation}
  \label{eq:sphere-area}
  \pi_q = 2 \pi^{q/2}/\Gamma(q/2).
  \end{equation}
A point $\bm x \in S_{q-1}$ can be written in the polar form about the
north pole $\bm e_1 = (1,0, \ldots, 0)^T$ as
\begin{equation}
    \label{eq:polar-x}
    \bm x = \pm \begin{bmatrix} \cos \phi \\ \sin \phi \, \bm u \end{bmatrix},
    \quad \ 0 \leq \phi \leq \pi,
  \end{equation}
  where $\bm u$ is a unit
  $(q-1)$-dimensional vector.  If $q=2$ then $u=\pm1$ is just a scalar.

  Using the polar representation (\ref{eq:polar-x}), the surface measure
  on $S_{q-1}$, written $[d\bm x]$, say, can be written recursively as
  \begin{equation}
    \label{eq:dx}
    [d\bm x] = \sin^{q-2}\phi\,  d \phi\, [d\bm u].
  \end{equation}
  When $q=2$, the formula simplifies to $[d\bm x] = d\phi$. However,
  note (\ref{eq:acg2}) used a slightly different convention for
  $\phi$; the scalar $u = \pm1$ was not present and the angle
  $\phi$ was allowed to range through the whole circle,
  $-\pi < \phi \leq \pi$.

  For all dimensions $q \geq 2$, changing $\phi$ to $\pi-\phi$ and
  $\bm u$ to $-\bm u$ changes $\bm x$ to $-\bm x$.  Hence when studying
  antipodally symmetric p.d.f.s, it is sufficient to restrict $\phi$
  to the range $0 \leq \phi < \pi/2$.

  Let $\bm y$ be another point in $S_{q-1}$ with polar representation
  \begin{equation}
    \label{eq:polar-y}
    \bm y = \begin{bmatrix}
      \cos \theta \\  \sin \theta\, \bm u \end{bmatrix}.
  \end{equation}
  If $\bm u$ is the same as in (\ref{eq:polar-x}) and $\theta = 2\phi$,
  then $\bm y$ can be said to be obtained from $\bm x$ by \emph{doubling
    the angle}, where ``angle'' here means the colatitude $\phi$.  In
  dimensions $q>2$ the concept of doubling the angle is less general
  than the squaring operation on the circle ($q=2$) given in
  (\ref{eq:squaring}). In particular, when $q>2$ the operation of
  doubling the angle depends on the choice of north pole.

  For use below, consider the following linear function of a
  $q$-dimensional unit vector $\bm y$,
  \begin{equation}
    \label{eq:P-linear}
    P(\bm y) = 1+\lambda^2  - 2 \lambda\, \bm y^T \bm \mu_0,
  \end{equation}
  and partition the unit vector $\bm \mu_0 = (\mu_1, \bm \mu_2^T )^T$
  in terms of a scalar and a $(q-1)$-vector.  Using (\ref{eq:polar-x}) and
  (\ref{eq:polar-y}), $P(\bm y)$ can be rewritten as a quadratic
  function of $\bm x$ as follows,
\begin{align}
  P(\bm y)
  &= 1+\lambda^2 -2\lambda \, \bm y^T \bm \mu_0 \notag \\
  &= (1+\lambda^2)  -2 \lambda \mu_1 \cos \theta
    - 2 \lambda (\bm \mu_2^T \bm u) \sin \theta \notag\\
  &= (1+\lambda^2)(\cos^2\phi + \sin^2\phi)  -
    2 \lambda \mu_1 (\cos^2\phi-\sin^2\phi)  -
    4 \lambda (\bm \mu_2^T \bm u) \sin \phi \cos \phi \notag\\
  &= (1+\lambda^2) (x_1^2+\bm x_2^T \bm x_2)
    -2 \lambda \mu_1 (x_1^2-\bm x_2^T \bm x_2)
    - 4 \lambda (\bm \mu_2^T \bm x_2) x_1 \notag \\
  &= (1+\lambda^2-2 \lambda \mu_1) x_1^2 +
    (1+\lambda^2+2 \lambda \mu_1) \bm x_2^T \bm x_2 -
    4\lambda (\bm \mu_2^T \bm x_2) x_1 \notag \\
  &= Q(\bm x), \text{  say,}  \label{eq:lin-to-quad}
\end{align}
a homogeneous quadratic form $\bm x^T \bm A \bm x$ with matrix
\begin{equation}
  \label{eq:pd}
  \bm A = \begin{bmatrix} 1+\lambda^2-2 \lambda \mu_1 & -2\lambda \bm \mu_2^T \\
  -2\lambda \bm \mu_2 & (1+\lambda^2+2 \lambda \mu_1) I_{q-1} \end{bmatrix}
\end{equation}
Since  $\bm \mu_0^T\bm \mu_0 = 1$, and
$|\lambda|<1$, it can be checked that  $\bm A$ is positive definite.

\section{Projections from the sphere to Euclidean space}
\label{sec:sphere-proj}
In this section we look at two standard tangent projections from the
sphere to the Euclidean space.  It is convenient to set up the
definitions and notation for all dimensions $q \geq 2$.  We can then
specialize to the case $q=2$ and describe how the projections are
connected to the transformations of Section \ref{sec:circle-trans}.

The first is \emph{gnomonic projection}, taking the open hemisphere
$H_{q-1} = \{\bm x \in S_{q-1} : x_1 > 0\}$ to $\mathbb{R}^{q-1}$.  If
$\bm x$ is a unit $q$-vector in the open hemisphere, it can be
written in the form (\ref{eq:polar-x}) where $0 \leq \phi < \pi/2$ and
$\bm u$ is a unit $(q-1)$-vector.  As in (\ref{eq:circ-rep}), let
$r = \tan \phi$.  Then the gnomonic projection is defined by
  \begin{equation}
  \label{eq:gnomonic}
  \bm v =  r \,\bm u = \frac{\sin \phi}{\cos \phi} \bm u =
  \frac{\sin \phi}{x_1} \bm u.
  \end{equation}

  The second is stereographic projection, taking the sphere $S_{q-1}$,
  minus the point at $-\bm e_1$, to $\mathbb{R}^{q-1}$.  If
  $\bm y \in S_{q-1}$ is a unit vector other than $-\bm e_1$, write it 
  in the form  (\ref{eq:polar-y}), where $-\pi < \theta < \pi$.  As in
  (\ref{eq:circ-s}), let $ s = \tan (\theta/2) $.  Then the
  \emph{stereographic projection} of $\bm y$ is defined by
  \begin{equation}
  \label{eq:stereo}
  \bm w = s \bm u = \frac{\sin(\theta/2)}{\cos(\theta/2)} \bm u =
  \frac{\sin \theta}{1+y_1} \,\bm u
\end{equation}
since $\sin \theta = 2 \sin (\theta/2) \cos (\theta/2)$ and
$1+y_1 = 1+\cos \theta = 2 \cos^2 (\theta/2)$.  

If $\bm y$ is obtained from $\bm x$ by angle doubling, then the two
projections are identical.  That is, if $\theta = 2\phi$, then $r=s$
and $\bm v = \bm w$.  However, the mapping of the uniform measure on
the sphere to Euclidean space is different for the two projections.
For gnomonic projection, the polar coordinate representation
$\bm v = r\, \bm u$ states that $r$ is the radial part of $\bm v$ so
that Lebesgue measure in the tangent space $\mathbb{R}^{q-1}$ is
related to the uniform measure on the sphere by
\begin{align}
  d\bm v &= r^{q-2} dr \,[d\bm u] \notag\\
       &= (\sin \phi/\cos \phi)^{q-2} (dr/d\phi)\  d\phi \,[du] \notag\\
  &= \cos^{-q}\phi \{\sin^{q-2} \phi\  d\phi \, [du]\} \notag\\
  &= \cos^{-q}\phi \,[d\bm x], \label{eq:cov-gnomonic}
\end{align}
using (\ref{eq:dx}) and $dr/d\phi = \sec^2\phi$.  On the other hand,
for stereographic projection, the polar coordinate representation
$\bm w = s \, \bm u$ implies
\begin{align}
  d\bm w &= s^{q-2} ds \,[d\bm u] \notag\\
         &= \{\sin(\theta/2)/\cos(\theta/2)\}^{q-2}\, (ds/d\theta)\,
           d\theta \, [du] \notag \\
         &= \frac{1}{2} \{\sin(\theta/2)/\cos(\theta/2)\}^{q-2}
           \{\cos (\theta/2)\}^{-2}
           \sin^{-(q-2)} \theta \{\sin^{q-2} \theta\, d\theta\,[d\bm u]\}
           \notag \\
         &= \left(\frac12\right)^{q-1} \cos^{-2(q-1)} (\theta/2) [d \bm y]
           \label{eq:cov-stereo}
\end{align}
since $ds/d\theta = (1/2) \sec^2(\theta/2)$ and $\sin \theta = 2 \sin (\theta/2) \cos(\theta/2)$.
Except on the circle $q=2$, the two differentials involve different
powers of $\cos(\theta/2) = \cos \phi$.

Since $d\bm v = d\bm w$ both represent Lebesgue measure in $\mathbb{R}^{q-1}$,
(\ref{eq:cov-gnomonic}) and (\ref{eq:cov-stereo}) can combined to
describe  effect of angle doubling on the sphere,
$$
[d\bm y] = 2^{q-1} \cos^{q-2}\phi\, [d \bm x].
$$
The reason for the cosine factor is straightforward to understand
intuitively.  For example, consider the case $q=3$ corresponding to
the usual sphere.  For a constant value of a colatitude, the the
longitude can range between $0$ and $2\pi$, and the corresponding
points on the sphere lie on a small circle.  If $\phi$ is near
$\pi/2$, the corresponding small circle for $\bm x$ is near the
equator, a circle with circumference $2\pi$.  However, the
corresponding value of $\theta = 2\phi$ is near $\pi$ and the
corresponding small circle for $\bm y$ lies near the south pole with
circumference close to 0.

\begin{figure}
\begin{center}
  \includegraphics[width=5in,height=3in]{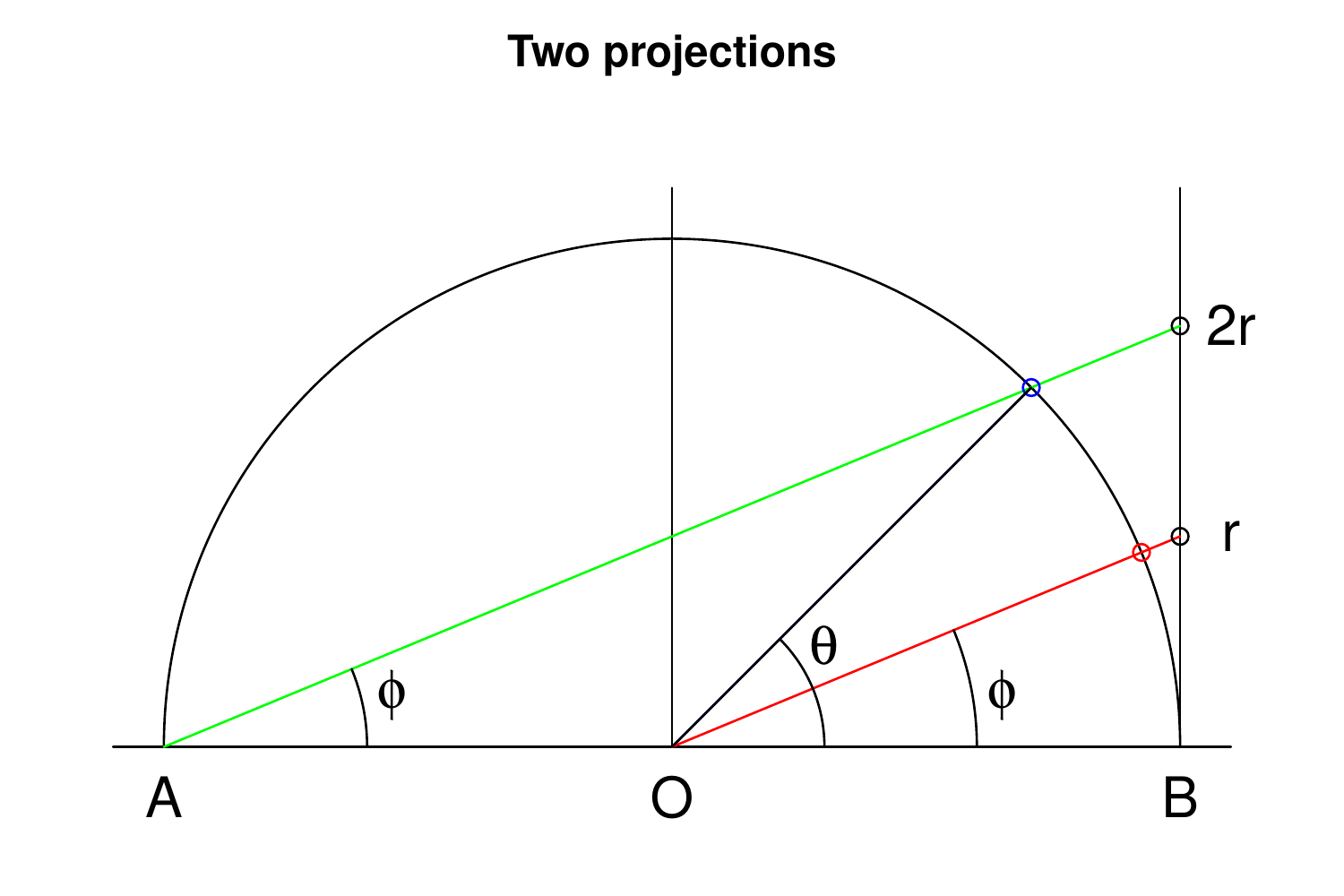}
  \end{center}
  \caption{Two projections, gnomonic and stereographic, from the
    circle to the vertical line tangent to the circle at point B.  If
    $\phi=\theta/2$, then $r = \tan \phi = \tan \theta/2$ is both the
    gnomonic projection of $\phi$ and the stereographic projection of
    $\theta$.}
    \label{fig:projs}
    
  \end{figure}

  Figure \ref{fig:projs} illustrates the two projections on the
  circle, where $\theta = 2\phi$.  The gnomonic projection of $\phi$
  is obtained by following the ray from the origin O through
  $(\cos \phi, \sin \phi)^T$ to the vertical line tangent to the
  circle at B.  Stereographic projection of $\theta$ is obtained by
  following the ray from A through $(\cos \theta, \sin \theta)^T$ to
  the same vertical line and dividing the result by 2.  Note the
  stereographic projection of $\theta$ is the same as the gnomonic
  projection of $\phi$.
  
  \section{The ACG distribution on the sphere}
  \label{sec:acg}
  This section takes a closer look at the ACG distribution on the sphere
  $S_{q-1}, \ q \geq 2$ and in particular derives its behavior under
  gnomonic projection.  First it is useful to recall some results
  about quadratic forms.
\subsection{Review of quadratic forms in the multivariate normal distribution}
\label{sec:acg:qf}
Let $\bm x = (\bm x_1^T, \bm x_2^T)^T$ be a $q$-dimensional vector
partitioned into two parts of dimensions $q_1$ and $q_2$.  Similarly
partition a $q \times q$ positive definite matrix as 
$$
\SSigma = \begin{bmatrix} \SSigma_{11} & \SSigma_{12} \\
  \SSigma_{21} & \SSigma_{22} \end{bmatrix}.
$$
If $\bm x$ follows a multivariate normal distribution,
$\bm x \sim N_q(\bm 0, \SSigma)$, then
$\bm x_1 \sim N_{q_1}(\bm 0, \SSigma_{11})$ and
$\bm x_2 | \bm x_1 \sim N_{q_2}(\bm \SSigma_{21} \SSigma_{11}^{-1} \bm
x_1, \SSigma_{22.1})$ \citep[e.g.][p. 63]{ka79}, where
$\SSigma_{22.1} = \SSigma_{22} - \SSigma_{21} \SSigma_{11}^{-1}
\SSigma_{12}$.  Writing the joint density of $\bm x$ as a product of a
marginal and a conditional density,
$f(\bm x) = f_1(\bm x_1) f(\bm x_2 | \bm x_1)$ yields an identity for
quadratic forms,
\begin{equation}
  \label{eq:qf-sum}
  Q = Q_1 + Q_{2.1}
  \end{equation}
where
\begin{align}
  Q &= \bm x^T \SSigma^{-1} \bm x \notag\\
  Q_1 &= \bm x_1^T \SSigma_{11}^{-1} \bm x_1, \label{eq:qf-details}\\
  Q_{2.1} &= (\bm x_2 - \SSigma_{21} \SSigma_{11}^{-1} \bm x_1)^T \SSigma_{22.1}^{-1}
(\bm x_2 - \SSigma_{21} \SSigma_{11}^{-1} \bm x_1) \notag.
\end{align}

If $q_1 = 1, \ q_2 = q-1$, then $\bm x_1 = x_1$ is a scalar,
$\SSigma_{11} = \bm \sigma_{11}$ is a scalar and
$\SSigma_{21} = \bm \sigma_{21}$ is a vector.  This case will be
useful in the next section when studying gnomonic projection.

  \subsection{Basic properties of the ACG distribution}
  \label{sec:acg:basic}
  This section reviews some basic facts about the ACG distribution.
Let $\SSigma$ be a symmetric $q \times q$ positive definite matrix
with inverse $\OOmega=\SSigma^{-1}$.  The angular central Gaussian (ACG)
distribution on $S_{q-1}$ is defined by the density (with respect to
the uniform measure on $S_{q-1}$) by
\begin{equation}
  \label{eq:acg}
  f_\ACG(\bm x) = f_\ACG(\bm x; \OOmega) =
  \pi_q^{-1} |\OOmega|^{1/2}/(\bm x^T \OOmega \bm x)^{q/2}.
\end{equation}
The parameter $\OOmega$ is defined up to a multiplicative scalar.  If
$\OOmega$ has spectral decomposition
$\OOmega = \bm \Gamma \bm \Delta \bm \Gamma^T$ where $\bm \Gamma$ is
an orthogonal containing the eigenvectors and $\bm \Delta $ is a
diagonal matrix containing the eigenvalues, then it is possible to
separate out the orientation and the concentration parts of the model.
The ACG distribution is antipodally symmetric,
$f_\ACG(\bm x) = f_\ACG(-\bm x)$.

If $q=2$ and $\OOmega =\diag(b^2,1)$ is a diagonal matrix with $0<b<1$,
then the density in polar coordinates reduces to (\ref{eq:acg2}).  A
similar expansion can be carried out in higher dimensions $q>2$.
Suppose $\OOmega$ is partitioned as 
$$
\OOmega = \begin{bmatrix} \omega_{11} &  \bm \omega_{21}^T \\ \bm \omega_{21}
  & \OOmega_{22} \end{bmatrix}
$$
and partition  a unit vector $\bm x \in S_{q-1}$ as in (\ref{eq:polar-x}).
The quadratic form becomes
\begin{equation}
  \label{eq:qf-acg1}
\bm x^T \OOmega \bm x = \omega_{11} \cos^2 \phi +
2 \sin \phi \cos \phi\, (\omega_{21}^T \bm u) +
\sin^2 \phi \, \bm u^T \OOmega_{22} \bm u.
\end{equation}
If, in addition, $\bm \omega_{21} = \bm 0$, then $\omega_{11}$ is an
eigenvalue.  If $\omega_{11}$ is the smallest eigenvalue, then the
density has its modes at $\phi = 0,\pi$.

\subsection{ACG distribution under gnomonic projection}
\label{sec:acg:proj}
Under gnomonic projection, (\ref{eq:qf-sum}) and
(\ref{eq:qf-details}) can be used to show that the ACG distribution on
the sphere is transformed to a multivariate Cauchy distribution in
$\mathbb{R}^{q-1}$.  To verify this result, recall the identities in
(\ref{eq:circ-rep}).  Then the quadratic form
$Q = Q(\bm x) = \bm x^T \SSigma^{-1} \bm x$, after dividing by
$\cos^2\phi = 1/(1+r^2)$, becomes
\begin{align}
  (1+r^2)Q &= \omega_{11} +  2 \bm v^T \bm \omega_{21} +
   \bm v^T \OOmega_{22} \bm v \notag\\
           &= \omega_{11} - \bm \omega_{21}^T \OOmega_{22}^{-1} \bm \omega_{21}
             + (\bm v + \OOmega_{22}^{-1} \bm \omega_{21})^T \OOmega_{22}
             (\bm v + \OOmega_{22}^{-1} \bm \omega_{21}) \notag \\
           &= \sigma_{11}^{-1} +
             (\bm v - \bm \sigma_{21}/ \sigma_{11})^T
             \SSigma_{22.1}^{-1} 
             (\bm v - \bm \sigma_{21} / \sigma_{11}), \label{eq:qf-acg}
\end{align}
using the identities
$\bm \sigma_{21}/\sigma_{11} = -\OOmega_{22}^{-1} \bm \omega_{21}$,
$\sigma_{11}^{-1} = \omega_{11} - \bm \omega_{21}^T \OOmega_{22}^{-1}
\bm \omega_{21}$ and $\SSigma_{22.1}^{-1} = \OOmega_{22}$ for the
inverse of a partitioned matrix \citep[e.g.,][p. 459]{ka79}.  Without
loss of generality we can rescale $\SSigma$ so that $\sigma_{11} = 1$.

The $(q-1)$-dimensional multivariate $t$-distribution, with location
parameter $\bm \mu$, scatter matrix $\bm B$ and degrees of freedom $\kappa>0$,
written $t_{q-1}(\bm \mu,\bm B,\kappa)$, has density proportional to 
\begin{equation}
  \label{eq:mvt}
  f(\bm v) \propto
  {\{1+\kappa^{-1} (\bm v - \bm \mu)^T \bm B^{-1}(\bm v - \bm \mu)\}^{-(q-1+\kappa)/2}}
\end{equation}
\citep[e.g.][]{ka79}.  If $\kappa=1$ the distribution is known as the multivariate
Cauchy distribution.

Using (\ref{eq:cov-gnomonic}), (\ref{eq:acg}) and (\ref{eq:qf-acg}) 
to give the p.d.f. of the ACG($\SSigma)$ 
  distribution after gnomonic projection yields 
  $$
  f_{\ACG,\text{gnomonic}}(\bm v) \propto
  Q^{-q/2} \cos^q \phi = Q^{-q/2} (1+r^2)^{-q/2},
  $$
  with respect to Lebesgue measure $d\bm v$ in the tangent plane, which is the same as (\ref{eq:mvt}) with $\kappa=1$.   That is,
  the gnomonic projection follows a multivariate Cauchy distribution
  $t_{q-1}(\bm \sigma_{21}, \SSigma_{22.1}^{-1},1)$.

\section{The spherical Cauchy distribution}
\label{sphere:sc}
\cite{Kato-McC20} have defined 
the \emph{spherical Cauchy (SC) distribution} on $S_{q-1}$ to have the p.d.f.
\begin{equation}
  \label{eq:sc}
  f_\text{SC}(\bm y; \lambda, \bm \mu_0) = \pi_q^{-1} \left\{\frac{1-\lambda^2}
    {P(\bm y)}\right\}^{q-1}, \quad P(\bm y) = 1+\lambda^2  - 2 \lambda \bm y^T \bm \mu_0,  \quad \bm y \in S_{q-1}.
  \end{equation}
  Here $0 \leq \lambda <1$ is a measure of concentration and
  $\bm \mu_0$ is a unit $q$-vector representing the modal direction.
  When $q=2$, the SC distribution reduces to the WC distribution 
  (\ref{eq:wc}).

  Write $\bm \mu_0 = (\mu_1, \bm \mu_2^T)^T$ where $\mu_1$ is a scalar
  and $\bm \mu_2$ is a $(q-1)$-vector and
  $\mu_1^2 + \bm \mu_2^T \bm \mu_2 = 1$.  Then, similarly to the expansion in 
  (\ref{eq:lin-to-quad}), the quantity $P(\bm y)$ in
  (\ref{eq:P-linear}) can be written in stereographic coordinates
  $\bm v$ as
\begin{align}
  P&=P(\bm y) = 1+\lambda^2 - 2 \lambda \bm y^T \bm \mu_0 \notag \\
  &=  (1+\lambda^2) - 2\lambda( \mu_1 \cos \theta + \bm u^T \bm \mu_2 \sin \theta) \notag \\
  &= \frac{1}{1+r^2}\{(1+\lambda^2)(1+r^2) -2\lambda [(1-r^2)\mu_1 + 2 \bm v^T \bm \mu_2)]\}\notag \\
  &= \frac{1}{1+r^2}\{\gamma + \delta r^2 - 4\lambda \bm v^T \bm \mu_2\}\notag \\
  &= \frac{1}{1+r^2}\{\gamma - (4 \lambda^2/\delta) \bm \mu_2^T \bm \mu_2 +
    \delta (\bm v - (2\lambda/\delta) \bm \mu_2)^T
           (\bm v - (2\lambda/\delta) \bm \mu_2)\}\notag \\
  &= \frac{\gamma^*}{1+r^2}\{1 + (\bm v - \bm m)^T (\bm v - \bm m)/\sigma^2\},
    \label{eq:sc-stereo}
\end{align}
where in the fourth line
$$
\gamma = 1+\lambda^2 - 2\lambda \mu_1, \quad
\delta = 1+\lambda^2 + 2\lambda \mu_1,
$$
and in the final line
$$
\gamma^* = \gamma - (4 \lambda^2/\delta) \bm \mu_2^T \bm \mu_2 =
(1-\lambda^2)/\delta, \quad 
\bm m = (2 \lambda/\delta) \,\bm \mu_2, \quad
\sigma = (1 - \lambda^2)/\delta.
$$
In addition the identities $\bm v^T\bm v = r^2 \bm u^T \bm u = r^2$,
$\cos^2 \phi= 1/(1-r^2)$,
$\cos \theta = (1-r^2)/(1+r^2)$, and
$\sin \theta = 2 \sin \phi \cos \phi = (2 \tan \phi)/ (1+r^2)$ have
been used.

Using the change of variables formula (\ref{eq:cov-stereo}), the
distribution of the stereographic projection of $\bm y$ has density
$$
f_{\gamma,\text{stereo}}(\bm v) \propto
P^{-(q-1)} \cos^{2(q-1)}(\theta/2) 
            = \{(1+r^2) P\}^{-(q-1)},
$$
which as a function of $\bm v$ can be identified with the density of
the multivariate $t$-distribution
$t_{q-1}(\bm m, (q-1)^{-1} \sigma^2 \bm I_{q-1},q-1)$ distribution with
$\kappa=q-1$ degrees of freedom.  Note the identification is valid even
$\bm \mu_0 \ne \bm e_1$, i.e. even if the the mode of the SC
distribution does not lie in the direction of the first coordinate
axis.  This result was proved in \citet{Kato-McC20}; see also
\citet{McC96} for a deeper study of the circular case.

Note the factor $(1+r^2)^{-(q-1)}$ in the density has the right power to combine with $P^{-(q-1)}$ in the density.  This property explains why the SC distribution was defined by raising
$P$ to the power $-(q-1)$, and not some other power, in (\ref{eq:sc}).

When $q \neq 2$, the SC distribution can never be identified
with the ACG distribution under angle doubling.  In particular, the
gnomonic projections of an ACG distribution follows a multivariate
Cauchy distribution (i.e. a multivariate $t$-distribution with 1
degree of freedom).  In contrast, the stereographic projection of an
SC distribution follows a multivariate $t$-distribution with $q-1$
degrees of freedom.

\section{Parameterizations and motivations for the
  wrapped Cauchy distribution on $S_1$}
\label{sec:wc-params}

The WC$(\lambda)$ distribution on the circle arises in a variety of
settings in statistics. Here we give a brief review.  The standard 
one-dimensional Cauchy distribution with scale parameter $b^2$ and
written $t_1(0, b^2,1)$ in (\ref{eq:mvt}), plays a key role in two of
the settings.

\begin{itemize}

\item[(a)] \emph{Angle doubling}.  This topic has been the main theme of the
  paper.  In particular, the WC$(\lambda)$ distribution can be
  obtained from the ACG$(b)$ distribution  by angle doubling, where
  $b$ and $\lambda$ are related by (\ref{eq:lambda-to-b}).
  
\item[(b)] \emph{Stereographic projection}.  As noted in Sections
  \ref{sec:acg}-\ref{sphere:sc}, the WC$(\lambda)$ distribution can be
  obtained from the Cauchy distribution by inverse stereographic
  projection when $b$ is related to $\lambda$ by
  (\ref{eq:lambda-to-b}).

\item[(c)] \emph{Wrapping}.  If $Z \sim t_1(0, b^2,1)$, set
  $\Theta = Z \text{ mod }2\pi$.  Recall the Cauchy distribution has
  Fourier transform
  $ \hat{f}(t) = \exp(-b |t|), \ t \in \mathbb{R}, $ and its
  wrapped version has Fourier coefficients
  $\hat{f}(m), \ m \in \mathbb{Z}$.  Since the WC$(\lambda)$
  distribution has Fourier coefficients,
  $\lambda^{|m|}, \ m \in \mathbb{Z}, $ it follows that
  $\Theta \sim \text{WC}(\lambda)$ distribution with
  $\lambda=\exp(-b)$.  Note this value of $\lambda$ is different from
  (b).

\item[(d)] \emph{AR(1) process}.  Consider the 
first-order autoregression AR(1) model  in time series,
$$
X_{t+1} = \lambda X_t + \epsilon_t, \quad t \in \mathbb{Z},
$$
where the innovation sequence $\{\epsilon_t\}$ consists of  independent identically distributed $N(0, \sigma^2_\epsilon)$ random
variables with $\epsilon_t$ independent of $X_s, \ s<t$.  For $|\lambda|<1$, the model describes a stationary Gaussian
process with spectral density (after standardizing it to be a probability
density) given by the WC$(\lambda)$ density.

\item[(e)] \emph{CAR(1) process}.  Consider the first-order conditional
  autoregression CAR(1) model, defined by the conditional
  distributions
$$
X_t | \{X_s, \ s \neq t\} \sim N(\alpha(X_{t-1}+X_{t+1}), \sigma^2_\eta),
$$
indexed by $t \in \mathbb{Z}$.  For $|\alpha|<1/2$, this model defines
a stationary process which is the same as the stationary AR(1)
process.  The parameters are related by
$\alpha=\lambda/(1+\lambda^2)$.
\end{itemize}

\begin{table}[t!]
  \begin{center}
    \caption{Various parameterizations of the wrapped Cauchy distribution}
    \label{table:wc}
    \begin{tabular}{cccccc}
      Number &   Parameter & $A$ & $B$ & $C$ & Setting \\ \hline
      1&   $0 \leq \lambda<1$ & $1-\lambda^2$ & $1+\lambda^2$ & $2 \lambda $
                                             & wrapped Cauchy, AR(1)\\
      2&   $0 <b \leq 1$ & $ 2b$ & $1+b^2$ & $1-b^2$ & doubled ACG, \\
      &&&&& stereographic projection\\
      3&   $0<\mu \leq \pi/2 $ & $\sin \mu$ & $1$ & $\cos \mu$ &  angular rep\\
      4&   $0 \leq  \alpha <1/2$ & $\sqrt{1-4\alpha^2}$ & $1$ & $2\alpha$ & CAR(1) \\
      \hline
    \end{tabular}
  \end{center} \end{table}

Several of these settings involve different ways to parameterize the
WC distribution.  Note that the WC$(\lambda)$ density for $0 \leq \lambda <1$
can be written in the form
\begin{equation}
  \label{eq:wc-ABC}
  f_\WC(\theta; \lambda) =  \frac{1}{2\pi} \frac{A}{B - C \cos \theta}, \quad
  \theta \in S_1,
\end{equation}
where $A,B>0$ and $C \geq 0$.  Provided $B^2 = A^2+C^2$, the density
integrates to 1.  Further, the density is unchanged if the parameters
are multiplied by the same scalar constant.  Hence, there is only one
free parameter.  Table \ref{table:wc} lists some common choices for
$A,B,C$. Further, by interchanging $A$ and $C$, as has already been
done for Parameterizations 1 and 2, the number of parameterizations
can be doubled.

Parameterization 1 is the standard representation.  As noted in (a),
Parameterization 2 is motivated by doubling the angle in the ACG
distribution with its standard parameterization.  As noted in (b), it
is also motivated by the standard parameterization of the Cauchy
distribution after inverse stereographic projection.  Parameterization
3 is the simplest algebraically.  Parameterization 4 is motivated the
the CAR(1) model in (e).

\bibliography{Kent}
\end{document}